\documentclass[11pt]{article}
\usepackage{palatino}
\usepackage{amsmath}
\usepackage{amsthm}
\usepackage{amssymb}
\usepackage{latexsym}
\usepackage{pstricks}
\usepackage{graphicx, pst-plot, pst-node, pst-text, pst-tree}
\usepackage{titlefoot}
\usepackage{titlesec}
\usepackage[small,it]{caption}
\usepackage{color}
\definecolor{lightgray}{rgb}{0.8, 0.8, 0.8}
\definecolor{darkgray}{rgb}{0.7, 0.7, 0.7}
\definecolor{darkblue}{rgb}{0, 0, .4}

\usepackage[bookmarks,breaklinks]{hyperref}
\hypersetup{
        colorlinks=true,
        linkcolor=darkblue,
        anchorcolor=darkblue,
        citecolor=darkblue,
        urlcolor=darkblue,
        pdfpagemode=UseThumbs,
        pdftitle={Subclasses of the Separable Permutations},
        pdfsubject={Combinatorics},
        pdfauthor={Albert, Atkinson, and Vatter},
}

\newtheorem{theorem}{Theorem}[section]
\newtheorem{proposition}[theorem]{Proposition}

\newtheorem{corollary}[theorem]{Corollary}

\newcounter{todocounter}

\newcommand{\Av}{\operatorname{Av}}

\newcommand{\Cl}{\operatorname{Cl}}

\newcommand{\C}{\mathcal{C}}
\newcommand{\D}{\mathcal{D}}
\newcommand{\E}{\mathcal{E}}

\renewcommand{\P}{\mathcal{P}}
\newcommand{\Q}{\mathcal{Q}}
\newcommand{\R}{\mathcal{R}}
\renewcommand{\S}{\mathcal{S}}
\newcommand{\Sep}{\mathcal{S}}
\newcommand{\T}{\mathcal{T}}
\newcommand{\U}{\mathcal{U}}

\newcommand{\X}{\mathcal{X}}
\newcommand{\Xinf}{\X}
\newcommand{\Xstraight}{\mathsf{X}}

\newcommand{\strongcomp}{\operatorname{sc}}

\newcommand{\vh}{{\mathbf h}}
\newcommand{\vbee}{{\mathbf b}}
\newcommand{\venn}{{\mathbf n}}
\newcommand{\vzero}{{\mathbf 0}}
\newcommand{\vone}{{\mathbf 1}}
\newcommand{\vpee}{{\mathbf p}}
\newcommand{\vque}{{\mathbf q}}
\newcommand{\vvee}{{\mathbf v}}

\makeatletter
\newfont{\footsc}{cmcsc10 at 8truept}
\newfont{\footbf}{cmbx10 at 8truept}
\newfont{\footrm}{cmr10 at 10truept}
\renewenvironment{abstract}%
                {
                  \begin{list}{}%
                     {\setlength{\rightmargin}{1in}%
                      \setlength{\leftmargin}{1in}}%
                   \item[]\ignorespaces\begin{small}}%
                 {\end{small}\unskip\end{list}}

\amssubj{05A05, 05A15}

\newpagestyle{main}[\small]{
        \headrule
        \sethead[\usepage][][]
        {\sc Subclasses of the Separable Permutations}{}{\usepage}}

\title{\sc Subclasses of the Separable Permutations}
\author{%
Michael H. Albert\\[-0.25ex]
\small Department of Computer Science\\[-0.5ex]
\small University of Otago\\[-0.5ex]
\small Dunedin, New Zealand\\[1.5ex]
M. D. Atkinson\\[-0.25ex]
\small Department of Computer Science\\[-0.5ex]
\small University of Otago\\[-0.5ex]
\small Dunedin, New Zealand\\[1.5ex]
Vincent Vatter\\[-0.25ex]
\small Department of Mathematics\\[-0.5ex]
\small Dartmouth College\\[-0.5ex]
\small Hanover, New Hampshire USA\\[1.5ex]
}

\titleformat{\section}
        {\large\sc}
        {\thesection.}{1em}{}   

\date{}

\begin{document}
\maketitle

\pagestyle{main}

\begin{abstract} 
  We prove that all subclasses of the separable permutations not
  containing $\Av(231)$ or a symmetry of this class have rational
  generating functions.  Our principal tools are partial well-order,
  atomicity, and the theory of strongly rational permutation classes
  introduced here for the first time.
\end{abstract}

\section{Introduction}
\label{sep-231-intro} 

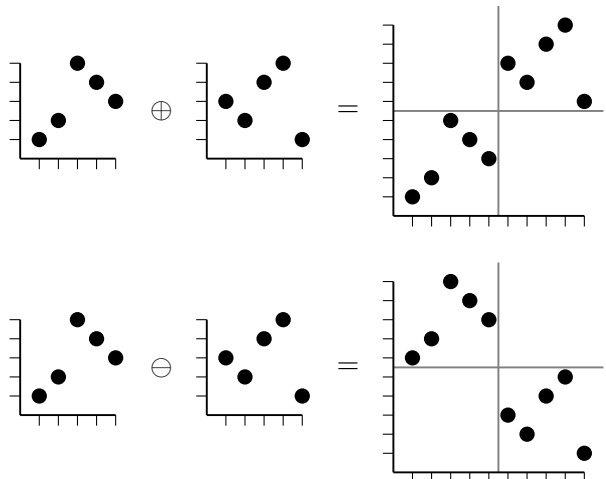
\begin{figure}
\[
\begin{array}{ccccc}
\psset{xunit=0.01in, yunit=0.01in}
\psset{linewidth=1\psxunit}
\begin{pspicture}(0,-30)(50,50) \psaxes[dy=10, Dy=1, dx=10, Dx=1,tickstyle=bottom, showorigin=false, labels=none](0,0)(50,50)
\pscircle*(10,10){4.0\psxunit}
\pscircle*(20,20){4.0\psxunit}
\pscircle*(30,50){4.0\psxunit}
\pscircle*(40,40){4.0\psxunit}
\pscircle*(50,30){4.0\psxunit}
\end{pspicture}
&
\psset{xunit=0.01in, yunit=0.01in}
\begin{pspicture}(0,0)(20,110)
\rput[c](10,55){$\oplus$}
\end{pspicture}
&
\psset{xunit=0.01in, yunit=0.01in}
\psset{linewidth=1\psxunit}
\begin{pspicture}(0,-30)(50,50) \psaxes[dy=10, Dy=1, dx=10, Dx=1,tickstyle=bottom, showorigin=false, labels=none](0,0)(50,50)
\pscircle*(10,30){4.0\psxunit}
\pscircle*(20,20){4.0\psxunit}
\pscircle*(30,40){4.0\psxunit}
\pscircle*(40,50){4.0\psxunit}
\pscircle*(50,10){4.0\psxunit}
\end{pspicture}
&
\psset{xunit=0.01in, yunit=0.01in}
\begin{pspicture}(0,0)(20,110)
\rput[c](10,55){$=$}
\end{pspicture}
&
\psset{xunit=0.01in, yunit=0.01in}
\psset{linewidth=1\psxunit}
\begin{pspicture}(0,0)(110,110)
\psline[linecolor=gray](55,0)(55,110)
\psline[linecolor=gray](0,55)(110,55)
\psaxes[dy=10, Dy=1, dx=10, Dx=1, tickstyle=bottom, showorigin=false, labels=none](0,0)(100,100)
\pscircle*(10,10){4.0\psxunit}
\pscircle*(20,20){4.0\psxunit}
\pscircle*(30,50){4.0\psxunit}
\pscircle*(40,40){4.0\psxunit}
\pscircle*(50,30){4.0\psxunit}
\pscircle*(60,80){4.0\psxunit}
\pscircle*(70,70){4.0\psxunit}
\pscircle*(80,90){4.0\psxunit}
\pscircle*(90,100){4.0\psxunit}
\pscircle*(100,60){4.0\psxunit}
\end{pspicture}
\\\\
\psset{xunit=0.01in, yunit=0.01in}
\psset{linewidth=1\psxunit}
\begin{pspicture}(0,-30)(50,50) \psaxes[dy=10, Dy=1, dx=10, Dx=1,tickstyle=bottom, showorigin=false, labels=none](0,0)(50,50)
\pscircle*(10,10){4.0\psxunit}
\pscircle*(20,20){4.0\psxunit}
\pscircle*(30,50){4.0\psxunit}
\pscircle*(40,40){4.0\psxunit}
\pscircle*(50,30){4.0\psxunit}
\end{pspicture}
&
\psset{xunit=0.01in, yunit=0.01in}
\begin{pspicture}(0,0)(20,110)
\rput[c](10,55){$\ominus$}
\end{pspicture}
&
\psset{xunit=0.01in, yunit=0.01in}
\psset{linewidth=1\psxunit}
\begin{pspicture}(0,-30)(50,50) \psaxes[dy=10, Dy=1, dx=10, Dx=1,tickstyle=bottom, showorigin=false, labels=none](0,0)(50,50)
\pscircle*(10,30){4.0\psxunit}
\pscircle*(20,20){4.0\psxunit}
\pscircle*(30,40){4.0\psxunit}
\pscircle*(40,50){4.0\psxunit}
\pscircle*(50,10){4.0\psxunit}
\end{pspicture}
&
\psset{xunit=0.01in, yunit=0.01in}
\begin{pspicture}(0,0)(20,110)
\rput[c](10,55){$=$}
\end{pspicture}
&
\psset{xunit=0.01in, yunit=0.01in}
\psset{linewidth=1\psxunit}
\begin{pspicture}(0,0)(110,110)
\psline[linecolor=gray](55,0)(55,110)
\psline[linecolor=gray](0,55)(110,55)
\psaxes[dy=10, Dy=1, dx=10, Dx=1, tickstyle=bottom, showorigin=false, labels=none](0,0)(100,100)
\pscircle*(10,60){4.0\psxunit}
\pscircle*(20,70){4.0\psxunit}
\pscircle*(30,100){4.0\psxunit}
\pscircle*(40,90){4.0\psxunit}
\pscircle*(50,80){4.0\psxunit}
\pscircle*(60,30){4.0\psxunit}
\pscircle*(70,20){4.0\psxunit}
\pscircle*(80,40){4.0\psxunit}
\pscircle*(90,50){4.0\psxunit}
\pscircle*(100,10){4.0\psxunit}
\end{pspicture}
\end{array}
\]
\caption{An example of a direct sum and a skew
sum.}\label{fig-ex-sum-skew}
\end{figure}

The \emph{separable permutations} are those which can be built from the permutation $1$
by repeatedly applying two operations, known as \emph{direct sum} (or
simply, sum) and \emph{skew sum} which are defined, respectively,
on permutations $\pi$ of length $m$ and $\sigma$ of length $n$ by
\[
\begin{array}{lcll} 
  (\pi\oplus\sigma)(i) &=&
  \left\{ 
    \begin{array}{l}
      \pi(i) \\
      \sigma(i-m)+m
    \end{array}
  \right.  &
  \begin{array}{l}
    \mbox{if $1\le i\le m$,} \\
    \mbox{if $m+1\le i \le m+n$,}
  \end{array} \\[15pt] 
  (\pi\ominus\sigma)(i) &=&
  \left\{
    \begin{array}{l}
      \pi(i)+n \\
      \sigma(i-m)
     \end{array}
  \right.  &
  \begin{array}{l}
    \mbox{if $1\le i\le m$,} \\
    \mbox{if $m+1 \le i \le m+n$,}
  \end{array}
\end{array}
\]
In this introductory section, we recapitulate some known results about
the separable permutations and some related sets of permutations.  The
operations $\oplus$ and $\ominus$ are best understood by considering
the plots of the permutations, as in Figure~\ref{fig-ex-sum-skew}.

While the term ``separable permutation'' dates only to the work of
Bose, Buss, and Lubiw~\cite{bose:pattern-matchin:}, these permutations first arose
in Avis and Newborn's work on pop stacks~\cite{avis:on-pop-stacks-i:}.
Separable permutations are the permutation analogues of two other
well-studied classes of object: complement-reducible graphs (also called cographs, or simply $P_4$-free graphs), and series-parallel (or $N$-free) posets.  A folkloric
result, which follows from the characterizations of these analogous
classes, characterizes the separable permutations.

\begin{proposition}
\label{prop-sep-basis} 
A permutation $\pi$ is separable if and only if it contains neither
$2413$ nor $3142$.
\end{proposition}

In Proposition~\ref{prop-sep-basis}, we say that the permutation $\pi$
of length $n$ \emph{contains} the permutation $\sigma$ of $[k]$
(written $\sigma\le\pi$) if $\pi$ has a subsequence of length $k$
order isomorphic to $\sigma$.  For example, $\pi=89167342$ (written
in list, or one-line notation) contains $\sigma=51342$, as can be seen
by considering the subsequence $91672$.

Our interest is with sets of permutations which are closed downward
under this containment order, which we call \emph{permutation classes}
(or occasionally just \emph{classes}).  Thus, $\C$ is a class if for
all $\pi$ in $\C$ and all $\sigma \le \pi$, $\sigma$ is also in $\C$.
One way to specify classes is as closures: if $X$ is any set of
permutations, its \emph{closure} is the permutation class
\[
\Cl(X)=\{\sigma : \mbox{$\sigma\le\pi$ for some $\pi\in X$}\}.
\]
However, it is often more useful to specify classes by what they do
not contain; for any permutation class $\C$ there is a unique (and
possibly infinite) antichain $B$ such that
\[
\C=\Av(B)=\{\pi: \pi \not \geq\beta\mbox{ for all } \beta \in B\}.
\]
This antichain $B$ is called the \emph{basis} of $\C$, so the basis of
the class of separable permutations is $\{2413,3142\}$.  Another way
to characterize the separable permutations is provided by the next
result; in this result we write $\C\oplus\D$ for the set of
permutations of the form $\pi\oplus\sigma$ where $\pi$ lies in $\C$
and $\sigma$ lies in $\D$, and extend this definition to $\C\ominus\D$
analogously.  (If $\C$ and $\D$ are permutation classes, then so are $\C \oplus \D$ and $\C \ominus \D$.)

\begin{proposition}\label{prop-sep-smallest} The class of separable
permutations is the smallest nonempty class $\C$ which satisfies both
$\C\oplus\C\subseteq\C$ and $\C\ominus\C\subseteq\C$.
\end{proposition}

Throughout this paper, $\Sep$ will be used to denote the class of
separable permutations.  Any class which satisfies
$\C\oplus\C\subseteq\C$ is called \emph{sum closed}, while any
class satisfying $\C\ominus\C\subseteq\C$ is called \emph{skew closed}.

The separable permutations contain a notable subclass, $\Av(231)$,
which Knuth~\cite{knuth:the-art-of-comp:1} showed are precisely the
permutations that can be sorted by a stack (a last-in first-out list).
A result similar to Proposition~\ref{prop-sep-smallest} holds for
$\Av(231)$ as well:

\begin{proposition}\label{prop-231-smallest} The class $\Av(231)$ is
the smallest nonempty class $\C$ which satisfies both $\C\oplus\C\subseteq\C$ and
$1\ominus\C\subseteq\C$.
\end{proposition}

In fact, since both basis elements ($2413$ and $3142$) of $\Sep$
each contain every non-monotone permutation of length $3$, all four of the
classes $\Av(132)$, $\Av(213)$, $\Av(231)$, and $\Av(312)$ are
contained in $\Sep$, and each of these has a characterisation similar
to the one given by Proposition~\ref{prop-231-smallest}. These four classes are all symmetric images of
one another under the operations of
reversal, inverse, and complementation (or compositions of these), each of
which preserves the containment order.

For any class $\C$ (or more generally any set of permutations), we denote by $\C_n$ the set of permutations in
$\C$ of length $n$, and say that the \emph{generating function} for
$\C$ is $\sum|\C_n|x^n$.  Whether this sum includes the empty
permutation ($n=0$) is a matter of taste and convenience, and we
generally elect not to include it.

Note that every separable permutation of length at least $2$ --- and
by extension, every permutation in $\Av(231)$ of length at least $2$
--- is either \emph{sum decomposable}, meaning that it is equal to
$\pi\oplus\sigma$ for two shorter (but nonempty) permutations $\pi$
and $\sigma$, or it is \emph{skew decomposable}, which is defined
analogously.  No permutation is both sum and skew decomposable, so the
separable permutations may therefore be partitioned into three sets:
$\{1\}$, the sum decomposable separable permutations, and the skew
decomposable separable permutations.  This observation allows one to
easily enumerate the class.

\begin{proposition}\label{prop-sep-gf} The generating function for the
separable permutations is
\[
\frac{1-x-\sqrt{1-6x+x^2}}{2}
\]
and thus the number of separable permutations of length $n$ is the
$n^{\mbox{\scriptsize th}}$ large Schr\"oder number.
\end{proposition}
\begin{proof} Let $f$ denote the generating function for the class of
separable permutations, $f_\oplus$ the generating function for its sum
decomposable elements, and $f_\ominus$ the generating function for its
skew decomposable elements.  As observed above, we have
$f=x+f_\oplus+f_\ominus$.  Any sum decomposable permutation may be
written uniquely as the direct sum of a sum indecomposable permutation
and another permutation, so since the class of separable permutations
is sum closed, we have $f_\oplus=(f-f_\oplus)f$.  Solving this shows
that $f_\oplus=f^2/(1+f)$, and then by symmetry $f_\ominus=f^2/(1+f)$,
so $f=x+2f^2/(1+f)$.  Solving this yields the desired generating
function.
\end{proof}

A similar approach gives the generating function for $\Av(231)$:

\begin{proposition}\label{prop-231-gf} The generating function for
$\Av(231)$ is
$$
\frac{1-2x-\sqrt{1-4x}}{2x}
$$
and thus the number of $231$-avoiding permutations of length $n$ is
the $n^{\mbox{\scriptsize th}}$ Catalan number.
\end{proposition}
\begin{proof} Let $f$ denote the generating function for $\Av(231)$,
and let $f_\oplus$ and $f_\ominus$ count the sum and skew decomposable
permutations in this class.  Since $\Av(231)$ is sum closed, by the
same logic as in the proof of Proposition~\ref{prop-sep-gf},
$f_\oplus=f^2/(1+f)$.  Now note that $\pi\ominus\sigma\in\Av(231)$ if
and only if $\pi$ is decreasing and $\sigma\in\Av(231)$.  Thus every
skew decomposable permutation in $\Av(231)$ may be written uniquely as
$1\ominus\sigma$ for $\sigma\in\Av(231)$, so $f_\ominus=xf$.
Substituting these values into our equation $f=x+f_\oplus+f_\ominus$
yields that $f=x(1+f)^2$, and solving this gives the generating
function claimed.
\end{proof}

Note that both of these generating functions are nonrational.
Clearly, we cannot hope for the generating function of a generic superclass
$\C\supseteq\Av(231)$ to be rational (although this may happen by
accident).  Our main result establishes the converse: if $\C$ is a
subclass of the separable permutations that does not contain any of
$\Av(132)$, $\Av(213)$, $\Av(231)$ or $\Av(312)$, then $\C$ has a
rational generating function.

\section{Partial Well-Order and Atomicity}\label{sec-pwo}

Many of our arguments depend on the \emph{partial well-order (pwo)} property.
In the context of the containment order on permutations, a permutation
class has the pwo property if it does not contain an infinite antichain.  This property has the following well-known consequence which is
important to us because it allows us to consider minimal counterexamples
within a pwo class.

\begin{proposition}\label{pwo-subclasses-dcc} The subclasses of a pwo
class $\C$ satisfy the {\it minimum condition\/}, i.e., every family 
of subclasses of a pwo class $\C$ has a minimal subclass under inclusion.
\end{proposition}

\begin{proof}
If there were a family of subclasses with no minimal element then we
could, inductively, find a strictly descending chain
\[\C^1\supsetneq
\C^2\supsetneq\cdots\]
of subclasses of $\C$.  For each $i\ge 1$, choose
$\beta_i\in\C^i\setminus\C^{i+1}$.  The set of minimal elements of
$\{\beta_1,\beta_2\ldots\}$ is an antichain and therefore finite.  Hence there
exists an integer $n$ such that $\{\beta_1,\beta_2\ldots,\beta_n\}$ contains these minimal elements.  In particular $\beta_m\leq \beta_{n+1}$ for some $m\leq n$, but $\C^{n+1}$ is a class, and therefore $\beta_{m}\in\C^{n+1}\subset \C^{m+1}$, a contradiction.
\end{proof}



%

For any class $\C$, we define its \emph{sum completion} $\bigoplus\C$ as the smallest sum closed class containing $\C$, and we define its \emph{strong completion}, $\strongcomp(C)$, as the smallest $\oplus$ and skew sum closed class containing $\C$.  In this notation \cite[Theorem 2.5]{atkinson:partially-well-:} states

\begin{proposition}\label{prop-sc-pwo} The sum completion and the strong completion of a pwo
class are pwo.
\end{proposition}

Since, clearly,  the separable permutations are the strong completion of the set $\{1\}$ we see that the separable class is pwo. 

Another key concept that we shall require is atomicity.  A permutation class is called \emph{atomic} if it cannot be written as the union of two proper subclasses.  The notion of atomicity was first studied (in a more general context) by Fra{\"{\i}}ss{\'e}~\cite{fraisse:sur-lextension-:}, who established several alternative characterizations of this property.  The only characterization we require features in our next proposition.  For a proof of this result in the context of permutations, we refer to Atkinson, Murphy and Ru\v{s}kuc~\cite[Theorem 1.2]{atkinson:pattern-avoidan:}.

\begin{proposition}\label{atomic-spine}
If $\C$ is an atomic class then there is a chain $\alpha_1\le\alpha_2\le\cdots$
of permutations in $\C$ such that $\C=\Cl(\{\alpha_1,\alpha_2,\dots\})$.
\end{proposition}

We refer to such a chain as a \emph{spine} for the class.

\section{Strongly Rational Classes}\label{sec-strong-rats}

Our main goal is to prove that if a subclass of the separable
permutations does not contain $\Av(231)$ or any of its symmetries,
then it and all of its subclasses have rational generating functions.
In this section we study this powerful property in its own right,
beginning by naming it: the permutation class $\C$ is \emph{strongly
rational} if it and all of its subclasses have rational generating
functions.  While strongly rational classes are naturally defined and
appear to be the ``correct'' context in which to state and prove the
tools of this section, they have received virtually no attention
before, and many conjectures remain.

\begin{proposition}\label{prop-strong-rat-union-inter} The union and
intersection of two strongly rational classes is strongly
rational.
\end{proposition}
\begin{proof} The intersection of two strongly rational classes is
contained in both of them, and so strongly rational by definition.
Now suppose that $\C$ and $\D$ are strongly rational and that
$\E\subseteq\C\cup\D$.  Since $\E=(\C\cap\E)\cup(\D\cap\E)$, we can
enumerate it by inclusion-exclusion; the generating function for $\E$
is the generating function for $\C\cap\E$ plus the generating function
for $\D\cap\E$ minus the generating function for $\C\cap\D\cap\E$.  As
$\C$ and $\D$ are strongly rational, all of these generating functions
are rational, so $\E$ has a rational generating function,
verifying that $\C\cup\D$ is strongly rational.
\end{proof}

We note that Proposition~\ref{prop-strong-rat-union-inter} does not
hold for classes with rational generating functions in general.
Neither does our next proposition, which follows from an argument of
Atkinson and Stitt~\cite{atkinson:restricted-perm:wreath} first
formalized by Murphy~{\cite[Chapter 9]{murphy:restricted-perm:}}
(although not in this context).

\begin{proposition}\label{prop-strong-rat-pwo}
Strongly rational classes are partially well ordered.
\end{proposition}
\begin{proof}
Suppose that the class $\C$ is not pwo.  Therefore it contains an infinite antichain, and in particular contains an infinite antichain $A\subseteq\C$ with at most one member of each length.  If $A_1\neq A_2$ are two subsets of $A$, then the two subclasses $\C\cap\Av(A_1)$ and $\C\cap\Av(A_2)$ have different enumerations.  To see this, suppose that $\alpha$ of length $k$ is the shortest permutation in one but not both of $A_1$ and $A_2$.  Then $\C\cap\Av(A_1)$ and $\C\cap\Av(A_2)$ contain the same permutations of length less than $k$ but differ by one in the number of permutations of length $k$.

  
Because $A$ is infinite, it follows that $\C$ has uncountably many subclasses with different generating functions.  These generating functions cannot all be rational, so $\C$ is not strongly rational.
\end{proof}

Our next step on the path to more powerful tools is the following.

\begin{proposition}

\label{prop-strong-rat-indecomps} 
If the class $\C$ is strongly rational, then each of the sets of its sum
indecomposable permutations, its sum decomposable permutations, its
skew indecomposable permutations and its skew decomposable
permutations have rational generating functions.
\end{proposition}

\begin{proof} 
  
It suffices to prove the claim for the sum indecomposable
permutations in $\C$ as the remaining cases follow by symmetry or
subtraction. If the claim were false then, because strongly rational
classes are pwo by Proposition~\ref{prop-strong-rat-pwo}, the
minimum condition of Proposition~\ref{pwo-subclasses-dcc}
shows that any counterexample would have a minimal subclass that was
also a counterexample.  Choose $\C$ to be such a minimal
counterexample.  By Proposition \ref{prop-sc-pwo} $\bigoplus\C$ is pwo and so the antichain of minimal  elements of the difference
$\left(\bigoplus\C\right)\setminus\C$ is a finite set, say $\{\beta_1,\dots,\beta_m\}$.  Clearly the
$\beta_i$ are nothing other than
the sum decomposable basis elements of $\C$.  
Suppose that
\[
\beta_i=\beta_{i,1}\oplus\beta_{i,2}\oplus\cdots\oplus\beta_{i,n_i}
\]
where the $\beta_{i,j}$s are sum indecomposable.  Now, for any
permutation $\pi$, let $\vbee(\pi)=(b_1,\dots,b_m)$ where for
each $i$, $\pi$ contains $\beta_{i,1}\oplus\cdots\oplus\beta_{i,b_i}$
but avoids
$\beta_{i,1}\oplus\cdots\oplus\beta_{i,b_i}\oplus\beta_{i,b_{i+1}}$.
Note that $\vbee(\pi)\le\venn-\vone=(n_1-1,\dots,n_m-1)$ for all
permutations $\pi\in\C$.  (Here and in what follows the order $\le$ on vectors
is the \emph{dominance order}, meaning that
$(p_1,\dots,p_m)\le(q_1,\dots,q_m)$ if and only if $p_i\le q_i$ for
all $1\le i\le m$.)

We now define a variety of generating functions:
\begin{itemize}
\item $f$ denotes the generating function for the class $\C$,
\item for a vector $\vpee$ of natural numbers, $f_{\vpee}$ denotes the
generating function for all permutations in $\C$ which avoid
$\beta_{i,p_i+1}\oplus\cdots\oplus\beta_{i,n_i}$ for all $i$,
\item $f_{\oplus}$ denotes the generating function for the sum
decomposable permutations in $\C$,
\item $f_{\not\oplus}$ denotes the generating function for the sum
indecomposable permutations in $\C$, and
\item for a vector $\vpee$ of natural numbers, $f_{\not\oplus}^\vpee$
denotes the generating function for the sum indecomposable
permutations in $\C$ with $\vbee(\pi)=\vpee$.
\end{itemize}

Note that $f_{\not\oplus}$, the generating function we wish to prove
rational, is the sum of the generating functions
$f_{\not\oplus}^\vpee$ for all $\vzero\le\vpee\le\venn-\vone$.  We now
claim that if
$\vzero\le\vpee<\venn-\vone$ then $f_{\not\oplus}^\vpee$ is
rational. We establish this claim by induction on the sum of the
entries of $\vpee$.
It is clearly true for the base case $f_{\not\oplus}^\vzero$, as
this function counts sum indecomposable elements of the proper
subclass $\{\pi\in\C : \vbee(\pi)=\vzero\}$ of $\C$.  For larger
$\vpee$, $f_{\not\oplus}^\vpee$ can be expressed as the difference
between the generating function for sum indecomposable elements of the
proper subclass $\{\pi\in\C : \vbee(\pi)\le\vpee\}$ of $\C$ (which is
rational by our choice of $\C$) and the sum of the generating
functions $f_{\not\oplus}^\vque$ for all $\vzero\le\vque<\vpee$ (which
are rational by induction).  This claim established, our goal is only
to show that $f_{\not\oplus}^{\venn-\vone}$ is rational.

Now we aim to enumerate $\C$, thereby expressing $f$, which is known
to be rational, in terms of the $f_{\not\oplus}^\vpee$ and $f_{\vpee}$
functions.  In the resulting equation $f_{\not\oplus}^{\venn-\vone}$
will be the only term not already known to be rational, yielding a
contradiction and completing the proof.

Consider the sum decomposable permutations of $\C$.  Each of these can
be expressed uniquely as $\theta\oplus\phi$ where $\theta$ is indecomposable and
so, counting how many of the summands in $\beta_i$ are contained in $\theta$, we see that $f_\oplus$ is the sum of $f_{\not\oplus}^\vpee
f_\vpee$ for all vectors $\vzero\le\vpee\le\venn-\vone$.  Of course,
$f_{\not{\oplus}}$ is the sum of the  generating functions $f_{\not\oplus}^\vpee$ and so we
obtain the equation
\[
f = f_{\not\oplus}+f_\oplus =
\displaystyle\sum_{\vzero\le\vpee\le\venn-\vone} f_{\not\oplus}^\vpee
+ \displaystyle\sum_{\vzero\le\vpee\le\venn-\vone}
f_{\not\oplus}^\vpee f_{\vpee} =
\displaystyle\sum_{\vzero\le\vpee\le\venn-\vone} f_{\not\oplus}^\vpee
(f_{\vpee}+1)
\]
Isolating the final term on the right hand side shows:
\[
f_{\not\oplus}^{\venn-\vone}({f_{\venn-\vone}+1}) =
f-\displaystyle\sum_{\vzero\le\vpee<\venn-\vone}
f_{\not\oplus}^\vpee (f_\vpee+1).
\]

As we have previously remarked, every generating function on the
right-hand side and in the second factor of the left hand side is rational, so $f_{\not\oplus}^{\venn-\vone}$ and
thus $f_{\not\oplus}$ must be rational as well.  This contradiction to
our choice of $\C$ completes the proof.
\end{proof}

It would be possible at this point to use
Proposition~\ref{prop-strong-rat-indecomps} to prove that if $\C$ is
strongly rational then $\bigoplus\C$ is strongly rational as well.
Instead, we provide a more powerful tool which we need for the main
theorem.

A permutation is said to be \emph{skew-merged} if it is the union of
an increasing subsequence and a decreasing subsequence.  The class of
skew-merged permutations was first studied by
Stankova~\cite{stankova:forbidden-subse:} in one of the earliest
papers on permutation patterns, and later enumerated by
Atkinson~\cite{atkinson:permutations-wh:}.  Stankova proved that the
skew-merged permutations have the basis $\{2143,3412\}$, a result
which can also be seen to follow from F\"oldes and Hammer's
characterization of split graphs~\cite{foldes:split-graphs:}.  Our
interest lies with the class of separable skew-merged permutations,
\[
\X=\Av(2143,2413,3142,3412).
\]

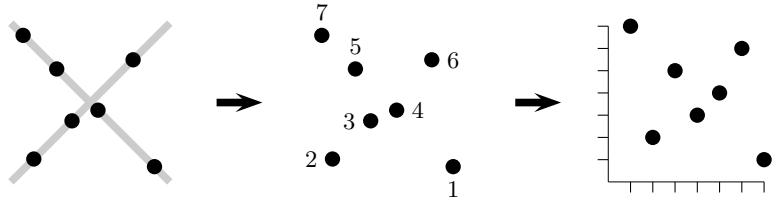
\begin{figure}
\begin{center}
\begin{footnotesize}
\begin{tabular}{ccccc}
\psset{xunit=0.008in, yunit=0.008in}
\psset{linewidth=0.005in}
\begin{pspicture}(0,0)(104,104)
\psline[linecolor=lightgray,linewidth=0.04in](0,0)(104,104)
\psline[linecolor=lightgray,linewidth=0.04in](0,104)(104,0)
\pscircle*(8,96){0.04in}
\pscircle*(15,15){0.04in}
\pscircle*(30,74){0.04in}
\pscircle*(40,40){0.04in}
\pscircle*(57,47){0.04in}
\pscircle*(80,80){0.04in}
\pscircle*(94,10){0.04in}
\end{pspicture}
&
\psset{xunit=0.008in, yunit=0.008in}
\psset{linewidth=0.04in}
\begin{pspicture}(0,0)(50,121)
\psline{->}(10,52)(40,52)
\end{pspicture}
&
\psset{xunit=0.008in, yunit=0.008in}
\psset{linewidth=0.005in}
\begin{pspicture}(0,0)(104,121)
\pscircle*(8,96){0.04in}
\pscircle*(15,15){0.04in}
\pscircle*(30,74){0.04in}
\pscircle*(40,40){0.04in}
\pscircle*(57,47){0.04in}
\pscircle*(80,80){0.04in}
\pscircle*(94,10){0.04in}
\rput[b](8,106){$7$}
\rput[r](5,15){$2$}
\rput[b](30,84){$5$}
\rput[r](30,40){$3$}
\rput[l](67,47){$4$}
\rput[l](90,80){$6$}
\rput[t](94,0){$1$}
\end{pspicture}
&
\psset{xunit=0.008in, yunit=0.008in}
\psset{linewidth=0.04in}
\begin{pspicture}(0,0)(50,121)
\psline{->}(10,52)(40,52)
\end{pspicture}
&
\psset{xunit=0.008in, yunit=0.008in}
\psset{linewidth=0.005in}
\begin{pspicture}(0,0)(104,104)
\psaxes[dy=14.571,Dy=1,dx=14.571,Dx=1,tickstyle=bottom,showorigin=false,  labels=none](0,0)(102,102)
\pscircle*(14.571,102){0.04in}
\pscircle*(29.143,29.143){0.04in}
\pscircle*(43.714,72.857){0.04in}
\pscircle*(58.286,43.714){0.04in}
\pscircle*(72.857,58.286){0.04in}
\pscircle*(87.429,87.429){0.04in}
\pscircle*(102,14.571){0.04in}
\end{pspicture}
\end{tabular}
\end{footnotesize}
\end{center}
\caption{The permutation $7253461$ can be drawn on an $\Xstraight$.}\label{fig-X-drawing}
\end{figure}

We label this class $\X$ because, in his thesis,
Waton~\cite{waton:on-permutation-:} showed that these are precisely
the permutations that can be ``drawn on an $\Xstraight$'' via the following
procedure: choose --- from an $\Xstraight$ made of right angles which form
$45^o$ angles with the axes in the plane --- $n$ points, no two lying
on a common vertical or horizontal line, and label these
points $1,\ldots,n$ reading bottom-to-top, then record these values
reading left-to-right, as depicted in Figure~\ref{fig-X-drawing}.
We may also define the class $\X$ in a manner similar to Propositions~\ref{prop-sep-smallest} and \ref{prop-231-smallest}.

\begin{proposition}\label{prop-X-smallest}
The class $\X$ is the smallest nonempty class $\C$ which contains $\C\oplus 1$, $1\oplus\C$, $\C\ominus 1$ and $1\ominus\C$.
\end{proposition}

Waton enumerated the class $\X$, obtaining the generating function
$(1-3x)/(1-4x+2x^2)$.  Later,
Elizalde~\cite{elizalde:the-x-class-and:} constructed a bijection
between the class $\X$ and the set of ``almost-increasing
permutations'' considered by Knuth~\cite[Section 5.4.8, Exercise
8]{knuth:the-art-of-comp:3}.

For the proof of our main result, we are interested not in the class $\X$ but rather in the $\X$-inflation of a strongly rational class $\U$.  This inflation, denoted $\X[\U]$, can be visualised by taking any permutation in $\X$, drawing it on the $\Xstraight$ as above, and then replacing each point in this drawing with a set of points corresponding to a permutation in $\U$ in such a way that the relationships between elements belonging to different points (of the permutation from $\X$) are the same as those between the original points.  Thus each point on the original drawing is ``inflated'' into a permutation from $\U$.  As we show below, such inflations are strongly rational.  Note, as made explicit in Corollary~\ref{cor-strong-rat-sum}, that this theorem implies that the sum completion of a strongly rational class is again strongly rational.

\begin{theorem}\label{thm-X-strong-rat} If $\U$ is a strongly rational
class then $\Xinf[\U]$ is also strongly rational.
\end{theorem}

\begin{proof} Let $\U$ be a strongly rational class.  It is
instructive to first consider the enumeration of $\Xinf[\U]$ itself.
Given a sum decomposable permutation in $\Xinf[\U]$, it may decompose in
one of two ways, either as a member of
$(\U_{\not\oplus})\oplus\Xinf[\U]$, or as a member of
$\Xinf[\U]\oplus(\U_{\not\oplus})$, or both, where $\U_{\not\oplus}$
denotes the set of sum indecomposable elements of $\U$.  The intersection of
these two sets is
$(\U_{\not\oplus})\oplus\{\Xinf[\U]\cup\epsilon\}\oplus(\U_{\not\oplus})$,
where $\epsilon$ denotes the empty permutation.  Doing the same for
skew decomposable elements of $\Xinf[\U]$ leads us to the equation
\begin{equation}\label{X-inflation-eqn}
g=x+2f_{\not\oplus}g-f_{\not\oplus}^2(g+1)+2f_{\not\ominus}g-f_{\not\ominus}^2(g+1),
\end{equation} where $g$ denotes the generating function for $\Xinf[\U]$,
$f_{\not\oplus}$ the generating function for $\U_{\not\oplus}$, and
$f_{\not\ominus}$ the generating function for $\U_{\not\ominus}$.
Solving for $g$ shows that it is indeed rational in $f_{\not\oplus}$
and $f_{\not\ominus}$, which are themselves elements of
$\mathbb{Q}(x)$ by
Proposition~\ref{prop-strong-rat-indecomps}. Specifically,
\[
g=\frac{x-f_{\not\oplus}^2-f_{\not\ominus}^2}{1-2f_{\not\oplus}+f_{\not\oplus}^2-2f_{\not\ominus}+f_{\not\ominus}^2}.
\]
Reassuringly, substituting $f_{\not\oplus}=f_{\not\ominus}=x$ gives us
the generating function, $(x-2x^2)/(1-4x+2x^2)$, which, upon adding
 $1$ to count the empty permutation, agrees with Waton's
enumeration of $\X=\Xinf[1]$.

In order to complete the proof --- that \emph{all} subclasses of $\Xinf[\U]$
have rational generating functions --- we adapt some notation of
Brignall, Huczynska, and Vatter~\cite{brignall:simple-permutat:}.  A
{\it property\/} is any set of permutations, and we say that $\pi$
{\it satisfies\/} the property $P$ if $\pi\in P$.  Given a set of
properties $\P$, we say that $\P$ is {\it separable query-complete\/}
if, for all nonempty permutations $\sigma$ and $\tau$ (not necessarily
lying in any class) and $P\in\P$, it can be decided whether
$\sigma\oplus\tau$ and $\sigma\ominus\tau$ satisfy $P$ given only the
knowledge about what properties in $\P$ are satisfied by $\sigma$ and
$\tau$.  For example, letting $\oplus$ denote the set of sum
decomposable permutations, we see that $\{\oplus\}$ is trivially
separable complete: assuming that $\sigma$ and $\tau$ are nonempty, $\sigma\oplus\tau$ always satisfies $\oplus$
while $\sigma\ominus\tau$ never satisfies $\oplus$.  Also note that
for any permutation $\beta$, the set $\{\Av(\delta) :
\delta\le\beta\}$ is separable query-complete: $\sigma\oplus\tau$ lies
in $\Av(\delta)$ if and only if $\sigma\in\Av(\gamma)$ or
$\tau\in\Av(\iota)$ for all $\gamma,\iota\le\delta\le\beta$ satisfying
$\gamma\oplus\iota=\delta$.

Returning to the situation at hand, consider an arbitrary subclass $\D\subseteq\Xinf[\U]$.
As $\U$ is strongly rational it is pwo by Proposition~\ref{prop-strong-rat-pwo}.  Thus
$\Xinf[\U]$ is contained in the strong completion of a pwo class, and so
is itself pwo by Proposition~\ref{prop-sc-pwo}.    Hence $\Xinf[\U]\setminus\D$ has only a finite number of minimal elements.  This set $B$ of minimal elements completely defines $\D$ as a subclass of  $\Xinf[\U]$ because $\Xinf[\U]\setminus\D$ is the set of permutations of $\Xinf[\U]$ that contain one or more permutations of $B$. It follows that
$\{\Av(\delta):\delta\in\Cl(B)\}$ is a \emph{finite} separable query-complete
set, since it is the union of a finite number of separable
query-complete properties.  Slightly more generally,
\[
\P=\{\oplus,\ominus\}\cup \{\Av(\delta):\delta\in\Cl(B)\}
\]
is also a finite separable query-complete set of properties.

For any permutation $\pi$, let $\P(\pi)$ denote the set of properties
in $\P$ satisfied by $\pi$.  We introduce three families of generating
functions which are defined for any subset $\Q\subseteq\P$:
\begin{itemize}
\item $f_\Q$, the generating function for the set $\{\pi\in\U :
\P(\pi)=\Q\}$,
\item $g_\Q$, the generating function for the set $\{\pi\in\Xinf[\U] :
\P(\pi)\supseteq\Q$\},
\item $h_\Q$, the generating function for the set $\{\pi\in\Xinf[\U] :
\P(\pi)=\Q\}$.
\end{itemize} Our goal, with this notation, is to show that all
functions of the form $g_\Q$ are rational, since it will then follow
that the generating function for $\D$, namely $g_{\{\Av(\beta):\beta\in
B\}}$, is rational.  First, note that all of the $f_\Q$ generating
functions are rational, by the strong rationality of $\U$ and
Proposition~\ref{prop-strong-rat-indecomps}.  Also note that $g_\Q$ is
the sum of all $h_\R$ with $\Q\subseteq\R\subseteq\P$.  Thus it will suffice to establish that the $h_\Q$ generating functions are
rational.

We are now ready to describe the analogues of the terms of
(\ref{X-inflation-eqn}) which relate the $h$ generating functions to
the $f$ generating functions and to each other.  Consider any subset
$\Q\subseteq\P$ of properties containing $\oplus$.  The permutations
satisfying $\Q$ must be sum decomposable, and thus can be expressed as
$\sigma\oplus\tau$ where at least one of $\sigma$ or $\tau$ is a sum
indecomposable permutation in $\U$.  Following our derivation of
(\ref{X-inflation-eqn}), such a $h_\Q$ can be expressed as a
linear combination of terms of three forms:
\begin{itemize}
\item $f_\R h_\S$ with $\oplus\notin\R$, which count permutations from
$(\U_{\not\oplus})\oplus \Xinf[\C]$,
\item $h_\S f_\T$ with $\oplus\notin\T$, which count permutations from
$\Xinf[\U]\oplus (\U_{\not\oplus})$, and
\item $f_\R h_\S f_\T$ with $\oplus\notin\R,\T$, which count
permutations from
$(\U_{\not\oplus})\oplus\Xinf[\U]\oplus(\U_{\not\oplus})$,
\end{itemize} the latter occurring with negative coefficients to
correct for over-counting.  Similarly, if $\Q\subseteq\P$ contains
$\ominus$, then $h_\Q$ can be expressed as a linear combination of
terms of the form $f_\R h_\S$ with $\ominus\notin\R$, $h_\S f_\T$ with
$\ominus\notin\T$, and $f_\R h_\S f_\T$ with $\ominus\notin\R,\T$.  As
no permutation can be both sum decomposable and sum indecomposable
there is only one more case, where neither $\oplus$ nor $\ominus$ lie
in $\Q$.  However, in this case the only permutations in $\Xinf[\U]$ that
can satisfy precisely the properties $\Q$ are those of $\C$, so here $h_\Q=f_\Q$.

Therefore, letting $\vh$ denote the column vector consisting of the
$h_\Q$ generating functions, there is some matrix $M$ of rational
functions in $\mathbb{Q}(x)$ and some constant vector $\vvee$ over
$\mathbb{Q}(x)$ such that $\vh=M\vh + \vvee$. Since all of our generating functions enumerate
non-empty permutations, the entries of $M$ all have zero constant term. Hence $I - M$ is invertible over $\mathbb{Q}(x)$, and thus each entry of $\vh$ is a rational function, proving the theorem.
\end{proof}

\begin{corollary}\label{cor-strong-rat-sum} If $\U$ is a strongly
rational class then $\bigoplus\U$ is also strongly rational, and if
$\D$ and $\E$ are strongly rational, then $\D\oplus\E$ is also
strongly rational.
\end{corollary}
\begin{proof} Since $\bigoplus\U\subseteq\Xinf[\U]$, the first part of
the corollary follows immediately from Theorem~\ref{thm-X-strong-rat}.
For the second part of the corollary,
Proposition~\ref{prop-strong-rat-union-inter} shows that $\D\cup\E$ is
strongly rational, and so $\bigoplus\left(\D\cup\E\right)$ is also
strongly rational.  However
$\D\oplus\E\subseteq\bigoplus\left(\D\cup\E\right)$ completing the proof.
\end{proof}

\section{Proof of the Main Result}

The machinery of Sections~\ref{sec-pwo} and
\ref{sec-strong-rats} developed, we are now ready to state and prove
our main result.

\begin{theorem}
\label{sep-231-main-thm} 
If  $\C$ is a subclass of the separable permutations that does not contain any of
$\Av(132)$, $\Av(213)$, $\Av(231)$ or $\Av(312)$ then $\C$ has a rational generating function.
\end{theorem}

\begin{proof} 
Suppose otherwise.  Because the class of separable permutations
is partially well ordered, its subclasses satisfy the minimum condition of
Proposition~\ref{pwo-subclasses-dcc}, and we can therefore choose
among all the counterexamples a minimal class $\C$.
We use two properties of $\C$ repeatedly:
\begin{itemize}
\item all proper subclasses of $\C$ have rational generating functions because $\C$ is a minimal counterexample, and thus
\item $\C$ is atomic because otherwise it would be the union of two strongly rational classes and hence strongly rational by Proposition~\ref{prop-strong-rat-union-inter}.
\end{itemize}
Our proof that $\C$ does not exist begins by ruling out some easy cases.  The easiest possibility to rule out is when $\C$ is either a sum or skew sum of two proper subclasses, which is eliminated by Corollary~\ref{cor-strong-rat-sum}.

Next we dispense with the case that $\C$ is sum closed (the case that $\C$ is skew closed is similar).  In this case we define $\C_{\not\oplus}$ as the set of sum 
indecomposable elements of $\C$.  It must be the case that
$\Cl(\C_{\not\oplus})=\C$ for otherwise it would be a proper subclass
of $\C$ and so strongly rational; but then, again by Corollary~\ref{cor-strong-rat-sum},
the sum closure of $\Cl(\C_{\not\oplus})$ would also be strongly rational and so its subclass
$\C$ would be strongly rational, a contradiction.  In the same way,
$\Cl(\C_{\not\ominus})=\C$.  In other words,
every permutation in $\C$ is contained in both a sum
indecomposable permutation and a skew sum indecomposable permutation
of $\C$.



We consider any spine $\alpha_1,\alpha_2,\dots$ of $\C$ (in the sense
of Proposition~\ref{atomic-spine}).
Clearly $1\oplus\alpha_1,1\oplus\alpha_2,\dots$ is also a spine for
$\C$ because we are assuming that $\C$ is sum closed, and since $\Cl(\C_{\not\oplus})=\C$, each of
these permutations is contained in a sum indecomposable element of
$\C$.  However, because $\C$ contains only separable permutations, the
sum indecomposable permutations in $\C$ (of length at least $2$) are
precisely the skew decomposable permutations.  The only way that
$1\oplus\alpha_i$ can be contained in a sum indecomposable element of
$\C$ is if it embeds completely into one of the skew components of
such a permutation.  Thus for all $i$ either $(1\oplus\alpha_i)\ominus
1$ or $1\ominus(1\oplus\alpha_i)$ lies in $\C$.  As one or the other
of these possibilities must occur infinitely often, we see that either
$\C=1\ominus\C$ or $\C=\C\ominus 1$.
Proposition~\ref{prop-231-smallest} now shows that $\C$ contains
$\Av(231)$ or a symmetry, $\Av(312)$, a contradiction.  Similarly, we
reach a contradiction if $\C$ is skew closed.  

For the remainder of the proof we may therefore take $\C$ to be
 neither a sum or skew sum of two proper subclasses nor to be sum closed or skew closed. 
 To complete the proof we shall find a proper subclass $\U\subsetneq\C$ for which
 $\C \subseteq \Xinf[\U]$.  This would indeed be a contradiction because, by the minimality
 of $\C$, $\U$ would be strongly rational and therefore, by Theorem \ref{thm-X-strong-rat},
 $\Xinf[\U]$ would also be strongly rational.

We now construct a finite collection of proper subclasses of $\C$ whose union will yield the desired $\U$.  To do this we shall rely on the following characterisation of subclasses of $\Xinf[\U]$, which follows trivially from Proposition~\ref{prop-X-smallest}.

\begin{proposition}\label{obs-subclass-of-inflation}
Given classes $\C$ and $\U$, we have $\C\subseteq\Xinf[\U]$ if and only if, for every $\pi\in\C$, one of the following holds (for nonempty $\gamma$ and $\tau$):
\begin{itemize}
\item $\pi \in \U$,
\item $\pi = \gamma \oplus \tau$ with $\gamma \in \U$ or $\tau
\in \U$, or
\item $\pi = \gamma \ominus \tau$ with $\gamma \in \U$ or $\tau \in
\U$.
\end{itemize}
\end{proposition}

With the aim of mimicking the structural decomposition provided by this proposition, we begin by defining
\[
\C_{SW}=\{\sigma\in\C : \sigma\oplus\C\subseteq\C\}.
\]
Note that $\C_{SW}$ is a proper subclass of $\C$ because $\C$ is not
sum closed.  In fact, $\C_{SW}$ is the maximum subclass $\P$ of $\C$
such that $\P \oplus \C \subseteq \C$. Similarly define $\C_{NW}$
maximal such that $\C_{NW}\ominus\C=\C$, $\C_{NE}$ maximal such that
$\C\oplus\C_{NE}=\C$, and $\C_{SE}$ maximal such that
$\C\ominus\C_{SE}=\C$.  As $\C$ is neither sum nor skew sum closed,
these are all proper subclasses of $\C$ (and may indeed be
empty). These are the first four classes that will be placed within
$\U$.

Consider any sum decomposable element $\pi$ of $\C$ and write
$\pi = \gamma \oplus \tau$ (in arbitrary fashion). If $\gamma \in
\C_{SW}$ or $\tau \in \C_{NE}$ then the conditions of
Proposition~\ref{obs-subclass-of-inflation} are already met. 
So suppose now that we have $\gamma \notin \C_{SW}$ and $\tau \notin
\C_{NE}$ with $\gamma \oplus \tau \in \C$. Define
\[
\E_\gamma=\{\sigma\in\C : \gamma\oplus\sigma\in\C\}.
\]
Clearly $\E_\gamma$ is a subclass of $\C$, and it is proper because
$\gamma\notin\C_{SW}$ and thus $\gamma\oplus\C\neq\C$.  Now define
\[
\D_\gamma=\{\sigma\in\C : \sigma\oplus\E_\gamma\subseteq\C\},
\]
where $\sigma\oplus\E_\gamma$ denotes $\{\sigma\oplus\iota : \iota\in\E_\gamma\}$.  Again, $\D_\gamma$ is a subclass of $\C$ and it is proper since $\tau
\in \E_{\gamma}$ and $\tau \notin \C_{NE}$. Therefore, since $\C$ is
not a sum of two proper subclasses, $\D_\gamma\oplus\E_\gamma$ is a proper
subclass of $\C$ for all permutations $\gamma$ of the type being considered.

While there may be infinitely many permutations $\gamma$ of this type
the number of \emph{distinct} classes $\D_\gamma$ is
finite.  To see this let $B$ denote the (finite) basis of $\C$ and
consider the sets $\Cl(\gamma)\cap\Cl(B)$, of which there are but a
finite number.  Suppose we have two of them which happen to be equal,
say $\Cl(\gamma)\cap\Cl(B)=\Cl(\overline{\gamma})\cap\Cl(B)$.  Now a
permutation $\sigma$ fails to lie in $\E_{\gamma}$ if and only if
$\gamma\oplus\sigma\not\in\C$.  But this happens if and only if
$\gamma\oplus\sigma$ contains some $\beta_1\oplus\beta_2\in B$ with
$\beta_1\leq\gamma$ and $\beta_2\leq\sigma$.  This means that
$\beta_1\in\Cl(\gamma)\cap\Cl(B)=\Cl(\overline{\gamma})\cap\Cl(B)$ and
so $\beta_1\leq\overline{\gamma}$.  In turn this implies that
$\overline{\gamma}\oplus\sigma\not\in\C$ and hence $\sigma$ fails to
lie in $\E_{\overline{\gamma}}$ also.  In other words
$\E_{\gamma}=\E_{\overline{\gamma}}$.  But then
$\D_\gamma=\D_{\overline{\gamma}}$ also.

Therefore --- in addition to $\C_{SW}$, $\C_{SE}$, $\C_{NW}$ and $\C_{NE}$ --- we include in $\U$ the finitely many classes $\D_\gamma\oplus\E_\gamma$ arising from
decompositions of this type.   By repeating an analogous argument for skew decompositions $\lambda \ominus \theta$ with $\lambda \notin \C_{NE}$ and $\theta \notin \C_{SW}$ we again find finitely many classes, and also include these in $\U$.  Since $\C$ is atomic, $\C$ cannot be equal to a finite union of proper subclasses, so so $\U\neq\C$.  This choice of $\U$ ensures that the conditions of Proposition~\ref{obs-subclass-of-inflation} are met, and thus that $\C \subseteq \Xinf[\U]$ establishing the desired contradiction, and proving the theorem.

\end{proof}

\section{Open problems}

The most obvious question is the converse to our main result: is
there a subclass of the separable permutations containing $\Av(231)$
which has a rational generating function? ÊIn fact, we are not aware of any
finitely based permutation class, separable or otherwise, which contains $\Av(231)$ and has a
rational generating function.  (Although we can non-constructively
prove that there are infinitely based classes satisfying these
conditions in a manner similar to the proof of Proposition~\ref{prop-strong-rat-pwo}.)

More generally, we are hopeful that the notion of strongly rational classes introduced herein will prove interesting and important in future studies of permutation classes.  One question, inspired by Theorem~\ref{thm-X-strong-rat}, would be: is there a natural characterisation of the classes $\C$ such that $\C[\U]$ is strongly rational for all strongly rational classes $\U$?  A positive answer to this question would lend hope to the possibility of a characterisation of the strongly rational classes themselves.

\bigskip

\bibliographystyle{acm}
\bibliography{../refs}

\def\cprime{$'$}
\begin{thebibliography}{10}

\bibitem{atkinson:permutations-wh:}
{\sc Atkinson, M.~D.}
\newblock Permutations which are the union of an increasing and a decreasing
  subsequence.
\newblock {\em Electron. J. Combin. 5\/} (1998), Research paper 6, 13 pp.

\bibitem{atkinson:partially-well-:}
{\sc Atkinson, M.~D., Murphy, M.~M., and Ru{\v{s}}kuc, N.}
\newblock Partially well-ordered closed sets of permutations.
\newblock {\em Order 19}, 2 (2002), 101--113.

\bibitem{atkinson:pattern-avoidan:}
{\sc Atkinson, M.~D., Murphy, M.~M., and Ru{\v{s}}kuc, N.}
\newblock Pattern avoidance classes and subpermutations.
\newblock {\em Electron. J. Combin. 12}, 1 (2005), Research paper 60, 18 pp.

\bibitem{atkinson:restricted-perm:wreath}
{\sc Atkinson, M.~D., and Stitt, T.}
\newblock Restricted permutations and the wreath product.
\newblock {\em Discrete Math. 259}, 1-3 (2002), 19--36.

\bibitem{avis:on-pop-stacks-i:}
{\sc Avis, D., and Newborn, M.}
\newblock On pop-stacks in series.
\newblock {\em Utilitas Math. 19\/} (1981), 129--140.

\bibitem{bose:pattern-matchin:}
{\sc Bose, P., Buss, J.~F., and Lubiw, A.}
\newblock Pattern matching for permutations.
\newblock {\em Inform. Process. Lett. 65}, 5 (1998), 277--283.

\bibitem{brignall:simple-permutat:}
{\sc Brignall, R., Huczynska, S., and Vatter, V.}
\newblock Simple permutations and algebraic generating functions.
\newblock {\em J. Combin. Theory Ser. A 115}, 3 (2008), 423--441.

\bibitem{elizalde:the-x-class-and:}
{\sc Elizalde, S.}
\newblock The {$\mathcal{X}$}-class and almost-increasing permutations.
\newblock {\it Ann. Comb.\/}, to appear.

\bibitem{foldes:split-graphs:}
{\sc F{\"o}ldes, S., and Hammer, P.~L.}
\newblock Split graphs.
\newblock In {\em Proceedings of the Eighth Southeastern Conference on
  Combinatorics, Graph Theory and Computing (Louisiana State Univ., Baton
  Rouge, La., 1977)\/} (Winnipeg, Man., 1977), Utilitas Math., pp.~311--315.
  Congressus Numerantium, No. XIX.

\bibitem{fraisse:sur-lextension-:}
{\sc Fra{\"{\i}}ss{\'e}, R.}
\newblock Sur l'extension aux relations de quelques propri\'et\'es des ordres.
\newblock {\em Ann. Sci. Ecole Norm. Sup. (3) 71\/} (1954), 363--388.

\bibitem{knuth:the-art-of-comp:1}
{\sc Knuth, D.~E.}
\newblock {\em The art of computer programming. {V}ol. 1: {F}undamental
  algorithms}.
\newblock Addison-Wesley Publishing Co., Reading, Mass., 1968.

\bibitem{knuth:the-art-of-comp:3}
{\sc Knuth, D.~E.}
\newblock {\em The art of computer programming. {V}olume 3}.
\newblock Addison-Wesley Publishing Co., Reading, Mass.-London-Don Mills, Ont.,
  1973.
\newblock Sorting and searching, Addison-Wesley Series in Computer Science and
  Information Processing.

\bibitem{murphy:restricted-perm:}
{\sc Murphy, M.~M.}
\newblock {\em Restricted Permutations, Antichains, Atomic Classes, and Stack
  Sorting}.
\newblock PhD thesis, Univ. of St Andrews, 2002.

\bibitem{stankova:forbidden-subse:}
{\sc Stankova, Z.~E.}
\newblock Forbidden subsequences.
\newblock {\em Discrete Math. 132}, 1-3 (1994), 291--316.

\bibitem{waton:on-permutation-:}
{\sc Waton, S.}
\newblock {\em On Permutation Classes Defined by Token Passing Networks,
  Gridding Matrices and Pictures: Three Flavours of Involvement}.
\newblock PhD thesis, Univ. of St Andrews, 2007.

\end{thebibliography}

\end{document}